\documentclass[11pt,a4paper]{amsart}
\usepackage{amsmath,amssymb,amsthm,enumerate,color}
\usepackage{graphicx}
\usepackage{enumerate}
\usepackage[all]{xy}
\usepackage[active]{srcltx}

\title
{A new generalization of Hermite's reciprocity law}

\author{Leandro Cagliero}
\address{CIEM-FaMAF, Universidad Nacional de C\'ordoba}
\email{cagliero@famaf.unc.edu.ar}

\author{Daniel Penazzi}
\address{CIEM-FaMAF, Universidad Nacional de C\'ordoba}
\email{penazzi@famaf.unc.edu.ar}

\date{}

\thanks{This work is supported in part by CONICET, Foncyt and SECYT-UNC grants.}

\numberwithin{equation}{section}

\theoremstyle{plain}
\newtheorem{theorem} {Theorem} [section]
\newtheorem{lemma} [theorem] {Lemma}

\theoremstyle{definition}

\newtheorem{remark} [theorem] {Remark}

\newtheorem{notation} [theorem] {Notation}

\newtheorem*{nonumbertheorem} {Theorem}

\renewcommand \parallel {/\kern-3pt/}

\newcommand \CC {\mathbb C}
\renewcommand \SS {\mathbb S}

\newcommand \pf{\begin{proof}}
 \newcommand \eop{\end{proof}}

\newcommand\dnt [7]{\langle#1,#2,#3\vert\vert #4,#5,#6,#7\rangle}

\newcommand\dntg [7]{\langle #1_{#2},\ldots,#1_{#3}\vert\vert #4_{#5},\ldots,#4_{#6},#4_{#7}\rangle}

\newcommand\vvv [2]{\langle#1\vert\vert #2\rangle}

\def\rro #1#2#3.{\rho_{#1,#2}(#3)}

\font\german=eufm10
\newcommand\gperm{\hbox{\german S}}
\newcommand \Mnvl [3]{\SS_{(#1^#2)}\left(\SS_{(#3)}(\CC^2)\right)}
\newcommand \Cnumbers{\mathbf c}
\newcommand \hook{\mathbf h}

\newcommand \lvec[1]{\stackrel{\leftarrow}{#1}}
\newcommand \rvec[1]{\stackrel{\rightarrow}{#1}}

\newcommand\Mcambio{ }

\newcommand\finMcambio{ }

\newcommand\fincambio{\color{black}}

\begin{document}

\maketitle


\begin{abstract}
  
Given a partition $\lambda$ of $n$, the {\it Schur  functor} 
$\SS_\lambda$
associates to any complex vector space $V$, a subspace 
$\SS_\lambda(V)$ of $V^{\otimes n}$.
Hermite's reciprocity law, in terms of the Schur functor, states that
$
 \SS_{(p)}\left(\SS_{(q)}(\CC^2)\right)\simeq \SS_{(q)}\left(\SS_{(p)}(\CC^2)\right).
$
We extend this identity to many other identities of the type 
$\SS_{\lambda}\left(\SS_{\delta}(\CC^2)\right)\simeq \SS_{\mu}\left(\SS_{\epsilon}(\CC^2)\right)$.
  
\end{abstract}

%
%

\section{Introduction}
\label{sec:intro}

Hermite's reciprocity law states that
\[
 {\rm Sym}^p\left({\rm Sym}^q(\CC^2)\right)\simeq 
{\rm Sym}^q\left({\rm Sym}^p(\CC^2)\right)
\] 
as 
$\text{GL}(2,\CC)$-modules, for any pair of non-negative integers $p$ and $q$,
(see e.g. \cite{FH}, Exercise 6.18). 
Here ${\rm Sym}^n(V)$ is the homogeneous component of degree $n$ 
in the symmetric algebra of $V$. 
This identity can also be stated in terms of the {\it Schur  functor}.
Recall that given any partition $\lambda$ of $n$, the Schur  functor 
$\SS_\lambda$
associates to any complex vector space $V$, a subspace 
(also known as the \emph{Weyl module})
$\SS_\lambda(V)$ of $V^{\otimes n}$ (see e.g. \S6.1 in \cite{FH}. 
We give some details in subsection \S \ref{subec.SchurFunctor}).
For instance, 
if $\lambda=(n)$ then $\SS_\lambda(V)\simeq {\rm Sym}^n(V)$,
and if $\lambda=(1^n)$ then $\SS_\lambda(V)\simeq \Lambda^n(V)$.

Thus, in terms of Schur functors,
Hermite's reciprocity law states that 
\[
 \SS_{(p)}\left(\SS_{(q)}(\CC^2)\right)\simeq \SS_{(q)}\left(\SS_{(p)}(\CC^2)\right).
\]
This reciprocity law has been extended to more general plethysms involving rectangle
partitions by L. Manivel in \cite{M}.
More precisely a proof of Hermite's reciprocity law can be obtained from the
Cayley-Silvester formula (\cite{CayleySilvester}); 
this formula was extended by M. Brion in \cite{B} and Manivel used it 
to prove the following extension of 
Hermite's reciprocity law, valid for all positive integers $n,k,d$:
\[
\begin{matrix} 
\Mnvl nk{d+k-1}&\simeq& \Mnvl dn{k+n-1}&\simeq& \Mnvl kd{n+d-1}\\[2mm]
\rule[-.36mm]{.6pt}{2.7mm}\,\wr && \rule[-.36mm]{.6pt}{2.7mm}\,\wr &&
\rule[-.36mm]{.6pt}{2.7mm}\,\wr\\[2mm]
\Mnvl nd{d+k-1}&\simeq& \Mnvl dk{k+n-1}&\simeq& \Mnvl kn{n+d-1}
\end{matrix}
\]
where the isomorphisms are now only as SL$(2,\CC)$-modules.

It is now natural to ask for other solutions to the following plethysm equation
\begin{equation}\label{eq.Main}
 \SS_{\lambda}\left(\SS_{\delta}(\CC^2)\right)\simeq \SS_{\mu}\left(\SS_{\epsilon}(\CC^2)\right)
\end{equation}
considering the partitions 
$\lambda$, $\delta$, $\mu$ and $\epsilon$ as unknowns
and the isomorphism either as SL$(2,\CC)$ or GL$(2,\CC)$-modules.

In this paper, we obtain new solutions to the 
plethysm equation \eqref{eq.Main} involving partitions of arbitrary 
number of `steps'. 
Manivel's result (involving rectangular partitions) turns out to be our 
one-step case.
In addition, we address the question of when an $\text{SL}(2,\CC)$-isomorphism 
is (or can twisted to obtain) an $\text{GL}(2,\CC)$-isomorphism.

\subsection*{Main results}
Let us denote 
$\SS_{\lambda}\left(\SS_{(d)}(\CC^2)\right)$ by $Y_{d+1}$ where 
$Y$ is the Young diagram of $\lambda$
(recall that 
$\dim \SS_{(d)}(\CC^2)=d+1$).
For instance
\begin{center}
 \setlength{\unitlength}{6pt}
\begin{picture}(11,5)(5,-4.3)
  \put(0,0){\line(1,0){3}}
  \put(0,-1){\line(1,0){3}}
  \put(0,-2){\line(1,0){2}}
  \put(0,-3){\line(1,0){2}}
  \put(0,-4){\line(1,0){1}}
  \put(0,0){\line(0,-1){4}}
  \put(1,0){\line(0,-1){4}}
  \put(2,0){\line(0,-1){3}}
  \put(3,0){\line(0,-1){1}}
  \put(2.3,-3.8){${}_z$}
  \put(4.5,-2){$=\,\SS_{\lambda}\left(\SS_{(z-1)}(\CC^2)\right),\quad\lambda=(3,2^2,1)$.}
\end{picture}
\end{center}
We add labels to a Young diagram to indicate the width and hight of the boxes.
For instance, if $\lambda=(9^2,5^4,3^4)$, its Young diagram is  

\begin{center}
\setlength{\unitlength}{9pt}
\begin{picture}(5,4.0)(-0,-3)
{\scriptsize
  \put(0,0){\line(1,0){3}}
  \put(0,-1){\line(1,0){3}}
  \put(0,-2){\line(1,0){2}}
  \put(0,-3){\line(1,0){1}}
  \put(0,0){\line(0,-1){3}}
  \put(1,0){\line(0,-1){3}}
  \put(2,0){\line(0,-1){2}}
  \put(3,0){\line(0,-1){1}}
  \put(0.2,.2){$3$}
  \put(1.2,.2){$2$}
  \put(2.2,.2){$4$}
  \put(-0.7,-0.7){$2$}
  \put(-0.7,-1.7){$4$}
  \put(-0.7,-2.7){$4$}
}  
\put(3.3,-1.7){.}
\end{picture}
\end{center}
One of the main results of the paper is the following theorem (see Theorem \ref{thm.MainMain}).
\begin{nonumbertheorem}
Let $x_1,\dots,x_n$ and $y_1,\dots,y_{n}$ be two sequences in $\mathbb{Z}_{\ge0}$, 
set $|x|=\sum x_i$, $|y|=\sum y_i$, and let $u,v,z\in\mathbb{Z}_{\ge0}$.
Then the following
$\text{SL}(2,\CC)$-isomorphism holds:

 \setlength{\unitlength}{12pt}
\begin{picture}(5,9.0)(-4,-8)
{\scriptsize
  \put(0,0){\line(1,0){7}}
  \put(0,-1){\line(1,0){7}}
  \put(0,-2){\line(1,0){6}}
  \put(0,-3){\line(1,0){5}}
  \put(0,-4){\line(1,0){4}}
  \put(0,-5){\line(1,0){3}}
  \put(0,-6){\line(1,0){2}}
  \put(0,-7){\line(1,0){1}}
  \put(0,0){\line(0,-1){7}}
  \put(1,0){\line(0,-1){7}}
  \put(2,0){\line(0,-1){6}}
  \put(3,0){\line(0,-1){5}}
  \put(4,0){\line(0,-1){4}}
  \put(5,0){\line(0,-1){3}}
  \put(6,0){\line(0,-1){2}}
  \put(7,0){\line(0,-1){1}}
  \put(0.2,.2){$x_1$}
  \put(1.2,.2){$\dots$}
  \put(2.2,.2){$x_n$}
  \put(3.2,.2){$\,u$}
  \put(4.2,.2){$y_1$}
  \put(5.2,.2){$\dots$}
  \put(6.2,.2){$y_n$}
  \put(-0.9,-0.7){$x_1$}
  \put(-0.9,-1.9){$\vdots$}
  \put(-0.9,-2.7){$x_n$}
  \put(-0.9,-3.7){$v$}
  \put(-0.9,-4.7){$y_1$}
  \put(-0.9,-5.9){$\vdots$}
  \put(-0.9,-6.7){$y_n$}
  \put(2.0,-6.9){$|x|+|y| + v +z$}
}
 \put(8,-2.7){$\simeq$}
\end{picture}
\begin{picture}(5,7.0)(-10,-8)
{\scriptsize
  \put(0,0){\line(1,0){7}}
  \put(0,-1){\line(1,0){7}}
  \put(0,-2){\line(1,0){6}}
  \put(0,-3){\line(1,0){5}}
  \put(0,-4){\line(1,0){4}}
  \put(0,-5){\line(1,0){3}}
  \put(0,-6){\line(1,0){2}}
  \put(0,-7){\line(1,0){1}}
  \put(0,0){\line(0,-1){7}}
  \put(1,0){\line(0,-1){7}}
  \put(2,0){\line(0,-1){6}}
  \put(3,0){\line(0,-1){5}}
  \put(4,0){\line(0,-1){4}}
  \put(5,0){\line(0,-1){3}}
  \put(6,0){\line(0,-1){2}}
  \put(7,0){\line(0,-1){1}}
  \put(0.2,.2){$x_1$}
  \put(1.2,.2){$\dots$}
  \put(2.2,.2){$x_n$}
  \put(3.2,.2){$\,v$}
  \put(4.2,.2){$y_1$}
  \put(5.2,.2){$\dots$}
  \put(6.2,.2){$y_n$}
  \put(-0.9,-0.7){$x_1$}
  \put(-0.9,-1.9){$\vdots$}
  \put(-0.9,-2.7){$x_n$}
  \put(-0.9,-3.7){$u$}
  \put(-0.9,-4.7){$y_1$}
  \put(-0.9,-5.9){$\vdots$}
  \put(-0.9,-6.7){$y_n$}
  \put(2.0,-6.9){$|x|+|y| + u + z$}
}
\end{picture}
\end{nonumbertheorem}

Although the diagrams in the above isomorphism have an odd number of steps,
it is immediate to derive from it (taking $u=0$) the following analogous 
isomorphism for even number of steps:

 \setlength{\unitlength}{12pt}
\begin{picture}(5,7.7)(-4,-6.7)
{\scriptsize
  \put(0,0){\line(1,0){6}}
  \put(0,-1){\line(1,0){6}}
  \put(0,-2){\line(1,0){5}}
  \put(0,-3){\line(1,0){4}}
  \put(0,-4){\line(1,0){3}}
  \put(0,-5){\line(1,0){2}}
  \put(0,-6){\line(1,0){1}}
  \put(0,0){\line(0,-1){6}}
  \put(1,0){\line(0,-1){6}}
  \put(2,0){\line(0,-1){5}}
  \put(3,0){\line(0,-1){4}}
  \put(4,0){\line(0,-1){3}}
  \put(5,0){\line(0,-1){2}}
  \put(6,0){\line(0,-1){1}}
  \put(0.2,.2){$x_1$}
  \put(1.2,.2){$\dots$}
  \put(2.2,.2){$x_n$}
  \put(3.2,.2){$\,u$}
  \put(4.1,.2){$y_1$}
  \put(4.7,.2){$.$}
  \put(4.9,.2){$.$}
  \put(5.1,.2){$.$}
  \put(5.4,.2){$y_{n\!-\!1}$}
  \put(-0.9,-0.7){$x_1$}
  \put(-0.9,-1.9){$\vdots$}
  \put(-0.9,-2.7){$x_n$}
  \put(-0.9,-3.7){$v$}
  \put(-0.9,-4.6){$y_1$}
  \put(-0.8,-5.0){$.$}
  \put(-0.8,-5.2){$.$}
  \put(-0.8,-5.4){$.$}
  \put(-1.2,-5.9){$y_{n\!-\!1}$}
  \put(2.0,-5.9){$|x|+|y| + v + z$}
}
 \put(7.5,-2.7){$\simeq$}
\end{picture}
\begin{picture}(5,7.7)(-9.5,-6.7)
{\scriptsize
  \put(0,0){\line(1,0){6}}
  \put(0,-1){\line(1,0){6}}
  \put(0,-2){\line(1,0){5}}
  \put(0,-3){\line(1,0){4}}
  \put(0,-4){\line(1,0){3}}
  \put(0,-5){\line(1,0){2}}
  \put(0,-6){\line(1,0){1}}
  \put(0,0){\line(0,-1){6}}
  \put(1,0){\line(0,-1){6}}
  \put(2,0){\line(0,-1){5}}
  \put(3,0){\line(0,-1){4}}
  \put(4,0){\line(0,-1){3}}
  \put(5,0){\line(0,-1){2}}
  \put(6,0){\line(0,-1){1}}
  \put(0.2,.2){$x_1$}
  \put(1.2,.2){$\dots$}
  \put(2.2,.2){$x_n$}
  \put(3.2,.2){$\,v$}
  \put(4.1,.2){$y_1$}
  \put(4.7,.2){$.$}
  \put(4.9,.2){$.$}
  \put(5.1,.2){$.$}
  \put(5.4,.2){$y_{n\!-\!1}$}
  \put(-0.9,-0.7){$x_1$}
  \put(-0.9,-1.9){$\vdots$}
  \put(-0.9,-2.7){$x_n$}
  \put(-0.9,-3.7){$u$}
  \put(-0.9,-4.6){$y_1$}
  \put(-0.8,-5.0){$.$}
  \put(-0.8,-5.2){$.$}
  \put(-0.8,-5.4){$.$}
  \put(-1.2,-5.9){$y_{n\!-\!1}$}
  \put(2.0,-5.9){$|x|+|y| + u + z$}
}
\end{picture}

Let's say that two pairs $(\lambda,d)$ and $(\mu,e)$ are equivalent if 
$ \SS_{\lambda}\left(\SS_{(d)}(\CC^2)\right)\simeq \SS_{\mu}\left(\SS_{(e)}(\CC^2)\right)$.
From the above isomorphism it is possible to obtain another isomorphism
by using the fact that an $\text{SL}(2,\CC)$-module is isomorphic to its dual module.
This, in general, yields an equivalence class of four different pairs $(\lambda,d)$.
If we additionally assume in the previous theorem that 
$x_i=v$ and $y_i=z$ for all $i=1,\dots,n$, then we can make use of its result twice, and obtain an equivalence class of six different pairs $(\lambda,d)$. 
In the odd case with $n=0$, this equivalence class of six pairs corresponds to Manivel's Theorem.

The even and odd cases with $n=1$ state that 

\begin{center}
\setlength{\unitlength}{10pt}
\begin{picture}(5,9)(-8,-7.5)
{\scriptsize
  \put(0,0){\line(1,0){2}}
  \put(0,-1){\line(1,0){2}}
  \put(0,-2){\line(1,0){1}}
  \put(0,0){\line(0,-1){2}}
  \put(1,0){\line(0,-1){2}}
  \put(2,0){\line(0,-1){1}}
  \put(0.2,.2){$v$}
  \put(1.2,.2){$z$}
  \put(-0.7,-0.7){$v$}
  \put(-0.7,-1.7){$u$}
  \put(1.3,-2.1){${u\!+\!v\!+\!z}$}
}  
  \put(4.4,-1.1){$\simeq$}
  \put(1.0,-3.6){\line(0,1){.7}}
  \put(1.2,-3.5){$\wr$}
\end{picture}
\begin{picture}(5,4.0)(-10,-7.5)
{\scriptsize
  \put(0,0){\line(1,0){2}}
  \put(0,-1){\line(1,0){2}}
  \put(0,-2){\line(1,0){1}}
  \put(0,0){\line(0,-1){2}}
  \put(1,0){\line(0,-1){2}}
  \put(2,0){\line(0,-1){1}}
  \put(0.2,.2){$v$}
  \put(1.2,.2){$u$}
  \put(-0.7,-0.7){$v$}
  \put(-0.7,-1.7){$z$}
  \put(1.3,-2.1){${\!v\!+\!2z}$}
}  
  \put(4.4,-1.1){$\simeq$}
  \put(1.0,-3.6){\line(0,1){.7}}
  \put(1.2,-3.5){$\wr$}
\end{picture}
\begin{picture}(5,4.0)(-12,-7.5)
{\scriptsize
  \put(0,0){\line(1,0){2}}
  \put(0,-1){\line(1,0){2}}
  \put(0,-2){\line(1,0){1}}
  \put(0,0){\line(0,-1){2}}
  \put(1,0){\line(0,-1){2}}
  \put(2,0){\line(0,-1){1}}
  \put(0.2,.2){$u$}
  \put(1.2,.2){$z$}
  \put(-0.7,-0.7){$v$}
  \put(-0.7,-1.7){$v$}
  \put(1.3,-2.1){${2v\!+\!z}$}
}  
  \put(1.0,-3.6){\line(0,1){.7}}
  \put(1.2,-3.5){$\wr$}
\end{picture}
\begin{picture}(5,4.0)(8,-2.5)
{\scriptsize
  \put(0,0){\line(1,0){2}}
  \put(0,-1){\line(1,0){2}}
  \put(0,-2){\line(1,0){1}}
  \put(0,0){\line(0,-1){2}}
  \put(1,0){\line(0,-1){2}}
  \put(2,0){\line(0,-1){1}}
  \put(0.2,.2){$z$}
  \put(1.2,.2){$v$}
  \put(-0.7,-0.7){$z$}
  \put(-0.7,-1.7){$u$}
  \put(1.3,-2.1){${u\!+\!v\!+\!z}$}
}  
  \put(4.4,-1.1){$\simeq$}
\end{picture}
\begin{picture}(5,4.0)(6,-2.5)
{\scriptsize
  \put(0,0){\line(1,0){2}}
  \put(0,-1){\line(1,0){2}}
  \put(0,-2){\line(1,0){1}}
  \put(0,0){\line(0,-1){2}}
  \put(1,0){\line(0,-1){2}}
  \put(2,0){\line(0,-1){1}}
  \put(0.2,.2){$u$}
  \put(1.2,.2){$v$}
  \put(-0.7,-0.7){$z$}
  \put(-0.7,-1.7){$z$}
  \put(1.3,-2.1){${\!v\!+\!2z}$}
}  
  \put(4.4,-1.1){$\simeq$}
\end{picture}
\begin{picture}(5,4.0)(4,-2.5)
{\scriptsize
  \put(0,0){\line(1,0){2}}
  \put(0,-1){\line(1,0){2}}
  \put(0,-2){\line(1,0){1}}
  \put(0,0){\line(0,-1){2}}
  \put(1,0){\line(0,-1){2}}
  \put(2,0){\line(0,-1){1}}
  \put(0.2,.2){$z$}
  \put(1.2,.2){$u$}
  \put(-0.7,-0.7){$z$}
  \put(-0.7,-1.7){$v$}
  \put(1.3,-2.1){${2v\!+\!z}$}
}  
\end{picture}
\end{center}
and 
\begin{center}
\setlength{\unitlength}{10pt}
\begin{picture}(7,10.0)(-6,-9)
{\scriptsize
  \put(0,0){\line(1,0){3}}
  \put(0,-1){\line(1,0){3}}
  \put(0,-2){\line(1,0){2}}
  \put(0,-3){\line(1,0){1}}
  \put(0,0){\line(0,-1){3}}
  \put(1,0){\line(0,-1){3}}
  \put(2,0){\line(0,-1){2}}
  \put(3,0){\line(0,-1){1}}
  \put(0.2,.2){$v$}
  \put(1.2,.2){$u$}
  \put(2.2,.2){$z$}
  \put(-0.7,-0.7){$v$}
  \put(-0.7,-1.7){$v$}
  \put(-0.7,-2.7){$z$}
  \put(1.4,-2.9){$2v\!+\!2z$}
}  
\put(4.5,-1.5){$\simeq$}
  \put(1.1,-4.4){\line(0,1){.7}}
  \put(1.3,-4.3){$\wr$}
\end{picture}
\begin{picture}(7,10.0)(-6,-9)
{\scriptsize
  \put(0,0){\line(1,0){3}}
  \put(0,-1){\line(1,0){3}}
  \put(0,-2){\line(1,0){2}}
  \put(0,-3){\line(1,0){1}}
  \put(0,0){\line(0,-1){3}}
  \put(1,0){\line(0,-1){3}}
  \put(2,0){\line(0,-1){2}}
  \put(3,0){\line(0,-1){1}}
  \put(0.2,.2){$v$}
  \put(1.2,.2){$v$}
  \put(2.2,.2){$z$}
  \put(-0.7,-0.7){$v$}
  \put(-0.7,-1.7){$u$}
  \put(-0.7,-2.7){$z$}
  \put(1.4,-2.9){$u\!+\!v\!+\!2z$}
}  
\put(4.5,-1.5){$\simeq$}
  \put(1.1,-4.4){\line(0,1){.7}}
  \put(1.3,-4.3){$\wr$}
\end{picture}
\begin{picture}(7,10.0)(-6,-9)
{\scriptsize
  \put(0,0){\line(1,0){3}}
  \put(0,-1){\line(1,0){3}}
  \put(0,-2){\line(1,0){2}}
  \put(0,-3){\line(1,0){1}}
  \put(0,0){\line(0,-1){3}}
  \put(1,0){\line(0,-1){3}}
  \put(2,0){\line(0,-1){2}}
  \put(3,0){\line(0,-1){1}}
  \put(0.2,.2){$v$}
  \put(1.2,.2){$z$}
  \put(2.2,.2){$z$}
  \put(-0.7,-0.7){$v$}
  \put(-0.7,-1.7){$u$}
  \put(-0.7,-2.7){$v$}
  \put(1.4,-2.9){$u\!+\!2v\!+\!z$}
}  
  \put(1.1,-4.4){\line(0,1){.7}}
  \put(1.3,-4.3){$\wr$}
\end{picture}
\begin{picture}(7,6.0)(16,-3)
{\scriptsize
  \put(0,0){\line(1,0){3}}
  \put(0,-1){\line(1,0){3}}
  \put(0,-2){\line(1,0){2}}
  \put(0,-3){\line(1,0){1}}
  \put(0,0){\line(0,-1){3}}
  \put(1,0){\line(0,-1){3}}
  \put(2,0){\line(0,-1){2}}
  \put(3,0){\line(0,-1){1}}
  \put(0.2,.2){$z$}
  \put(1.2,.2){$u$}
  \put(2.2,.2){$v$}
  \put(-0.7,-0.7){$z$}
  \put(-0.7,-1.7){$z$}
  \put(-0.7,-2.7){$v$}
  \put(1.4,-2.9){$2v\!+\!2z$}
}  
\put(4.5,-1.5){$\simeq$}
\end{picture}
\begin{picture}(7,6.0)(-6,-9)
{\scriptsize
  \put(0,0){\line(1,0){3}}
  \put(0,-1){\line(1,0){3}}
  \put(0,-2){\line(1,0){2}}
  \put(0,-3){\line(1,0){1}}
  \put(0,0){\line(0,-1){3}}
  \put(1,0){\line(0,-1){3}}
  \put(2,0){\line(0,-1){2}}
  \put(3,0){\line(0,-1){1}}
  \put(0.2,.2){$z$}
  \put(1.2,.2){$v$}
  \put(2.2,.2){$v$}
  \put(-0.7,-0.7){$z$}
  \put(-0.7,-1.7){$z$}
  \put(-0.7,-2.7){$u$}
  \put(1.4,-2.9){$u\!+\!v\!+\!2z$}
}  
\put(4.5,-1.5){$\simeq$}
\end{picture}
\begin{picture}(7,6.0)(-6,-9)
{\scriptsize
  \put(0,0){\line(1,0){3}}
  \put(0,-1){\line(1,0){3}}
  \put(0,-2){\line(1,0){2}}
  \put(0,-3){\line(1,0){1}}
  \put(0,0){\line(0,-1){3}}
  \put(1,0){\line(0,-1){3}}
  \put(2,0){\line(0,-1){2}}
  \put(3,0){\line(0,-1){1}}
  \put(0.2,.2){$z$}
  \put(1.2,.2){$z$}
  \put(2.2,.2){$v$}
  \put(-0.7,-0.7){$z$}
  \put(-0.7,-1.7){$v$}
  \put(-0.7,-2.7){$u$}
  \put(1.4,-2.9){$u\!+\!2v\!+\!z$}
}  
\end{picture}
\end{center}
\vskip -2cm
These, and other corollaries, are obtained in \S\ref{Corollaries}.

Recall that given a partition $\lambda$ and a number $d\ge0$ the \emph{hook length} 
of $\lambda$ and the $d$-\emph{content} of $\lambda$ are, respectively, 
the following polynomials
\[
\hook_\lambda(q)=\prod[h(u)]_q, \qquad
\Cnumbers^d_\lambda(q)=\prod[d+1+c(u)]_q,
\]
where $[a]_q$ is the $q$-analog of $a$,
$h$ and $c$ are, respectively, the hook and the content functions
and both products run over the entries of the Young diagram of $\lambda$.
It is known (see e.g. \cite[Ch. 7]{S}) that the SL$(2,\CC)$-character of 
$\SS_{\lambda}\left(\SS_{\delta}(\CC^2)\right)$ is, up to 
a power of $q$, equal to 
\[
 P_\lambda^d(q)=\frac{\Cnumbers^d_\lambda(q)}{\hook_\lambda(q)}
\]
where $d=\delta_1-\delta_2$.

The following theorem translates the plethysm equation \eqref{eq.Main} 
in terms of $P$.
Although the results stated in this theorem might be known, 
we did not find an explicit reference to it, thus we prove it in 
\S3 (see Theorem \ref{Thm.SLGL}).
If $\lambda$ is a partition, then $|\lambda|$ denotes the sum of its parts.
\begin{nonumbertheorem}
 Let  $\delta=(\delta_1,\delta_2)$, $\epsilon=(\epsilon_1,\epsilon_2)$
and $d=\delta_1-\delta_2$, $e=\epsilon_1-\epsilon_2$.
Let $\lambda$, $\mu$ be partitions with $\ell(\lambda)\le d+1$ and $\ell(\mu)\le e+1$.
Then 
\medskip
\begin{enumerate}[(1)]
 \item $\SS_\lambda\left(\SS_\delta(\CC^2)\right)\simeq 
\SS_\mu\left(\SS_\epsilon(\CC^2)\right)$ as  $\text{SL}(2,\CC)$-modules if and only if 
\begin{equation*}
P^d_\lambda=P^e_\mu
\end{equation*}
and in this case $|\lambda|d-|\mu|e$ is even.
\medskip
\item  $\SS_\lambda\left(\SS_\delta(\CC^2)\right)\simeq 
\SS_\mu\left(\SS_\epsilon(\CC^2)\right)$ as  $\text{GL}(2,\CC)$-modules
if and only if, in addition to $P^d_\lambda=P^e_\mu$, 
it also holds
\begin{equation*}
|\delta||\lambda|=|\epsilon||\mu|.
\end{equation*}
\end{enumerate}
\end{nonumbertheorem}

%

\section{Technical background}
\label{sec.techback}

\subsection{Partitions}\label{partitions}
A \emph{partition} $\lambda$ of $n$  is an ordered sequence of nonnegative
integers $\lambda_1\ge\lambda_2\ge...$ with $|\lambda|=n$, where $|\lambda|=\sum\lambda_i$.
The $\lambda_i$'s are called the parts of the partition 
and the length $\ell(\lambda)$ of $\lambda$ is the number of non zero parts.
If $k\ge\ell(\lambda)$ then $\lambda$ will be denoted as
$\lambda=(\lambda_1,\dots,\lambda_{k})$
or by indicating multiplicities with exponential notation, for instance
$(4,4,3,1,1,1)=(4^2,3,1^3)$.
If $\lambda$ and $\mu$ are two partitions, we denote by $\lambda+\mu$
the partition whose parts are $(\lambda+\mu)_i=\lambda_i+\mu_i$.

To each partition $\lambda=(\lambda_1,\dots,\lambda_{k})$ of $n$ we associate
its \emph{Young diagram} $Y(\lambda)$ and its \emph{standard tableau} $T(\lambda)$:
$Y(\lambda)$ is the graphical arrangement consisting of $k$ left-justified rows of boxes, 
with $\lambda_i$ boxes in the $i$-th row, and 
$T(\lambda)$ is the assigment of the integers $1,2,\dots,n$ to the $n$ boxes 
of $Y(\lambda)$ obtained by writing the numbers $1,2,\dots,n$ starting 
on the first row and increasing to the right and then continuing on the second row, etc.
For example, if $\lambda=(3,2,2,1)$ then
\begin{center}
 \setlength{\unitlength}{12pt}
\begin{picture}(5,3.8)(0,-3.7)
  \put(-3,-1.5){$Y(\lambda)=$}
  \put(0,0){\line(1,0){3}}
  \put(0,-1){\line(1,0){3}}
  \put(0,-2){\line(1,0){2}}
  \put(0,-3){\line(1,0){2}}
  \put(0,-4){\line(1,0){1}}
  \put(0,0){\line(0,-1){4}}
  \put(1,0){\line(0,-1){4}}
  \put(2,0){\line(0,-1){3}}
  \put(3,0){\line(0,-1){1}}
\end{picture}
\begin{picture}(5,3.8)(-5,-3.7)
  \put(-3,-1.5){$T(\lambda)=$}
  \put(0,0){\line(1,0){3}}
  \put(0,-1){\line(1,0){3}}
  \put(0,-2){\line(1,0){2}}
  \put(0,-3){\line(1,0){2}}
  \put(0,-4){\line(1,0){1}}
  \put(0,0){\line(0,-1){4}}
  \put(1,0){\line(0,-1){4}}
  \put(2,0){\line(0,-1){3}}
  \put(3,0){\line(0,-1){1}}
  \put(0.3,-.9){1}
  \put(1.3,-.9){2}
  \put(2.3,-.9){3}
  \put(0.3,-1.9){4}
  \put(1.3,-1.9){5}
  \put(0.3,-2.9){6}
  \put(1.3,-2.9){7}
  \put(0.3,-3.9){8}
\end{picture}
\end{center}

The {\it transpose} of a partition is the partition $\lambda^t$ whose Young diagram is the transpose of the Young diagram
of $\lambda$.
For example the transpose of the partition $(3,2,2,1)$ is the partition $(4,3,1)$ as can be seen by the drawing above.

\subsection{Schur functor}\label{subec.SchurFunctor}
If $\lambda$ is a partition of $n$, 
two subgroups of the symmetric group $\gperm_n$ are associated to $T(\lambda)$:
\begin{align*}
 P_\lambda &=\{\sigma\in \gperm_n:\sigma\text{ preserves each row of $T(\lambda)$}\}, \\
Q_\lambda&=\{\sigma\in \gperm_n:\sigma\text{ preserves each column of $T(\lambda)$}\}.
\end{align*}
Following \cite{FH} we denote by $a_\lambda$, $b_\lambda$, $c_\lambda$ 
the following elements of the group algebra $\CC[\gperm_n]$:
\[
a_\lambda=\sum_{\sigma\in P_\lambda}\sigma,\qquad
b_\lambda=\sum_{\sigma\in Q_\lambda}\text{sgn}(\sigma)\sigma,\qquad
c_\lambda=a_\lambda b_\lambda.
\]
The element $c_\lambda$ is called the \emph{Young symmetrizer} associated to $\lambda$.
The permutation group $\gperm_n$ 
acts naturally on $V^{\otimes n}$ by 
$\sigma.(v_1\otimes...\otimes v_n)=v_{\sigma(1)}\otimes...\otimes v_{\sigma(n)}$. 
This action is naturally extended to an action of its group algebra $\CC[\gperm_n]$. 
The image of $V^{\otimes n}$ under the action of $c_\lambda$ is denoted $\SS_\lambda(V)$
and the map $V\mapsto \SS_\lambda(V)$ is called the \emph{Schur functor}. 
  
For instance:
\begin{itemize}
\item $\SS_{(n)}(V)\simeq {\rm Sym}^n(V)$,
\item $\SS_{(1^n)}(V)\simeq \Lambda^n(V)$
\item $\SS_\lambda(V)=0$ if $\lambda$ has more than $\dim(V)$ parts.
\end{itemize}

\subsection{Schur polynomials}\label{subsec.SchurPoly}
If $\lambda$ is a partition of $n$ and $k\ge\ell(\lambda)$, 
the \emph{Schur polynomial}  in $k$ variables  associated to $\lambda$
is 
\[
 s_\lambda(x_1,\dots,x_k)=\frac{\det(x_j^{\lambda_i+k-i})}{\det(x_j^{k-i})},
\]
This is a symmetric polynomial in $k$ variables of degree $n$ for any $k\ge\ell(\lambda)$.

The Schur polynomial has an interesting property that will be useful later:
given a partition $\lambda$ and $k\ge\ell(\lambda)$ we will denote by $\lambda'$ the partition whose Young diagram  is the complement of $Y(\lambda)$ in the
$(k\times\lambda_1)$-rectangle.
(This definition depends on $k$, though this fact  is not indicated in the notation).
That is,
\[
\lambda'=(\lambda_1-\lambda_k,\dots,\lambda_{1}-\lambda_{2})
\]
For example, for $k=6$ we have that
\begin{center}
 \setlength{\unitlength}{9pt}
\begin{picture}(5,5.3)(0,-5.4)
  \put(-5.5,-1.5){if $\;Y(\lambda)=$}
  \put(0,0){\line(1,0){3}}
  \put(0,-1){\line(1,0){3}}
  \put(0,-2){\line(1,0){2}}
  \put(0,-3){\line(1,0){2}}
  \put(0,-4){\line(1,0){1}}
  \put(0,0){\line(0,-1){4}}
  \put(1,0){\line(0,-1){4}}
  \put(2,0){\line(0,-1){3}}
  \put(3,0){\line(0,-1){1}}
\end{picture}
\begin{picture}(6,5.5)(-6,-5.4)
  \put(-7,-1.5){then $\;Y(\lambda')=$}
  \put(0,0){\line(1,0){3}}
  \put(0,-1){\line(1,0){3}}
  \put(0,-2){\line(1,0){3}}
  \put(0,-3){\line(1,0){2}}
  \put(0,-4){\line(1,0){1}}
  \put(0,-5){\line(1,0){1}}
  \put(0,0){\line(0,-1){5}}
  \put(1,0){\line(0,-1){5}}
  \put(2,0){\line(0,-1){3}}
  \put(3,0){\line(0,-1){2}}
\end{picture}
\end{center}
It is not difficult to prove (see Exercise 7.41 of \cite{S}) that 
\begin{equation}\label{eq.SchurDual} 
  (x_1\dots x_k)^{\lambda_1}s_\lambda(x_1^{-1},\dots,x_k^{-1})=s_{ \lambda' }(x_1,\dots,x_k).
\end{equation}

\subsection{   Polynomial    representations of GL$(V)$ and SL$(V)$}\label{subsec.PolyRep}
Let $V$ be a finite dimensional complex vector space of dimension $k$. 
A \emph{polynomial representation} of GL$(V)$ is a finite dimensional
representation of GL$(V)$ such that the
matrix entries (associated to a given basis) are given by polynomial functions on $V$.
It is well known that every polynomial representation of GL$(V)$ 
can be decomposed into irreducible subrepresentations. In particular,
$\SS_\lambda(V)$ is an irreducible GL$(V)$-subrepresentation of $V^{\otimes n}$
for all partitions $\lambda$ of $n$.
The \emph{highest weight theorem} states that $\lambda\mapsto \SS_\lambda(V)$ establishes a
one-to-one correspondence 
between the set of equivalence classes of irreducible polynomial representations of GL$(V)$
and the set of partitions $\lambda$ with $\ell(\lambda)\le k$,
see for instance \S6 in \cite{FH}.

Moreover    $\lambda\mapsto \SS_\lambda(V)$ also establishes a one-to-one correspondence 
between the set of equivalence classes of irreducible polynomial representations of SL$(V)$
and the set of partitions $\lambda$ with $\ell(\lambda)\le k-1$.
This follows from the following fact:
if 
\[
 \tilde\lambda= \lambda-(\lambda_k^k)= (\lambda_1-\lambda_k,\dots,\lambda_{k-1}-\lambda_{k})
\]
then $\SS_\lambda(V)\simeq\SS_{(\lambda_{k}^k)}(V)\otimes\SS_{\tilde\lambda}(V)$
as GL$(V)$-modules. But since $\SS_{(r^{k})}(V)$ is the 1-dimensional GL$(V)$-module 
corresponding to $\det^{r}$, then we obtain that $\SS_\lambda(V)\simeq\SS_{\tilde\lambda}(V)$ as SL$(V)$-modules. Note that $\tilde\lambda$ now has at most $k-1$ parts.

\subsection{Characters of  GL$(V)$-modules}
If $\pi$ is a polynomial representation of  GL$(V)$, the \emph{character} of 
$\pi$ is the function $\chi_\pi:\text{GL}(V)\to\CC$ defined by 
$\chi_\pi(g)=\text{tr}(\pi(g))$.
If $\pi_1$ and $\pi_2$ are two polynomial representations of 
GL$(V)$ then $\pi_1\simeq\pi_2$ if and only if they have the same character.
Similarly, $\pi_1\simeq\pi_2$ as SL$(V)$-modules 
if and only if $\chi_{\pi_1}|_{\text{SL}(V)}=\chi_{\pi_2}|_{\text{SL}(V)}$.
If $g\in\text{GL}(V)$ has eigenvalues
$\theta_1,\dots,\theta_{k}$ (counted with multiplicities), then it is known that
\begin{equation}\label{eq.SchurCaracter}
 \chi_{_{\SS_\lambda(V)}}(g)=s_\lambda(\theta_1,\dots,\theta_{k})
\end{equation}
for any partition $\lambda$ with $\ell(\lambda)\le k$.
(See for instance \cite{FH}).

Let $\delta=(\delta_1,\delta_2)$ be a partition with at most two parts
and let $d=\delta_1-\delta_2$.
We know that 
$\dim\SS_\delta\left(\CC^2\right)=d+1$ and, 
as a representation of $\text{SL}(2,\CC)$,
$\SS_{\delta}(\CC^2)$ corresponds to the irreducible representation of the Lie algebra
$\mathfrak{sl}(2,\CC)$ of highest weight $d$.

If $g\in\text{GL}(2,\CC)$ has eigenvalues $x_1$ and $x_2$ then it
follows from \eqref{eq.SchurCaracter} that
\[
\chi_{_{\SS_\delta\left(\CC^2\right)}}(g)=s_\delta(x_1,x_2)=
x_1^{\delta_1}x_2^{\delta_2}+
x_1^{\delta_1-1}x_2^{\delta_2+1}+\dots+
x_1^{\delta_2}x_2^{\delta_1}.
\]
Hence the eigenvalues of $g$ in 
${\SS_\delta\left(\CC^2\right)}$ are 
$\{x_1^{\delta_1}x_2^{\delta_2},\dots, x_1^{\delta_2}x_2^{\delta_1}\}$
(all with multiplicity 1) and thus,
if $\lambda$ is a partition with $\ell(\lambda)\le d+1$, 
then the character of 
$\SS_{\lambda}\left(\SS_{\delta}(\CC^2)\right)$ is the plethysm
\begin{equation}\label{eq.Plethysm}
  \chi_{_{\SS_{\lambda}\left(\SS_{\delta}(\CC^2)\right)}}(g)
  =
  s_\lambda(x_1^{\delta_1}x_2^{\delta_2},
x_1^{\delta_1-1}x_2^{\delta_2+1},\dots,
x_1^{\delta_2}x_2^{\delta_1}). 
\end{equation}
In particular, if  $g\in\text{SL}(2,\CC)$ with eigenvalues $x_1$ and $x_1^{-1}$ then
\begin{equation*}
  \chi_{_{\SS_{\lambda}\left(\SS_{\delta}(\CC^2)\right)}}(g) 
  =
  \chi_{_{\SS_{\lambda}\left(\SS_{(d)}(\CC^2)\right)}}(g)
  =
  s_\lambda(x_1^{d},x_1^{d-2},\dots,x_1^{-d}). 
\end{equation*}
This identity and \eqref{eq.SchurDual} imply 
that if $\lambda'$ is as in \S\ref{subsec.SchurPoly}
\Mcambio (with $k=d+1$) \finMcambio
then  
$\chi_{_{\SS_{\lambda}\left(\SS_{(d)}(\CC^2)\right)}}$
and 
$\chi_{_{\SS_{\lambda'}\left(\SS_{(d)}(\CC^2)\right)}}$
coincide in ${\text{SL}(2,\CC)}$ and therefore
we obtain:
\begin{theorem}\label{eq.SL-Dual}
 $$ \SS_{\lambda}\left(\SS_{(d)}(\CC^2)\right)\simeq\SS_{\lambda'}\left(\SS_{(d)}(\CC^2)\right)$$
as ${\text{SL}(2,\CC)}$-modules.
\end{theorem}
 This corresponds to the fact that 
$\SS_{\lambda'}\left(\SS_{(d)}(\CC^2)\right)$ and
$\SS_{\lambda}\left(\SS_{(d)}(\CC^2)\right)$ are dual to each other as 
${\text{SL}(2,\CC)}$-modules and every polynomial representation of 
${\text{SL}(2,\CC)}$ is isomorphic to its dual.

\subsection{The Hook-content formula}

%
%

Given a natural number $a$, let 
\[
 [a]=[a]_q=\frac{1-q^a}{1-q}=1+q+\dots+q^{a-1}
\]
be the $q$-analog of $a$. 
If $u=(i,j)$ is a box of the Young diagram of $\lambda$ let $c(u)=j-i$
and let $h(u)$ be the number of boxes directly below or directly to the 
right of $u$, including $u$ once. 
For example, we indicate in the following diagrams the values of $c$ and $h$ 
respectively:
\begin{center}
 \setlength{\unitlength}{13pt}
\begin{picture}(5,4)(0,-3.9)
 \put(0,0){\line(1,0){3}}
  \put(0,-1){\line(1,0){3}}
  \put(0,-2){\line(1,0){2}}
  \put(0,-3){\line(1,0){2}}
  \put(0,-4){\line(1,0){1}}
  \put(0,0){\line(0,-1){4}}
  \put(1,0){\line(0,-1){4}}
  \put(2,0){\line(0,-1){3}}
  \put(3,0){\line(0,-1){1}}
  \put(0.3,-.9){\small \,0}
  \put(1.3,-.9){\small \,1}
  \put(2.3,-.9){\small \,2}
  \put(0.3,-1.9){\small \!-1}
  \put(1.3,-1.9){\small \,0}
  \put(0.3,-2.9){\small \!-2}
  \put(1.3,-2.9){\small \!-1}
  \put(0.3,-3.9){\small \!-3}
\end{picture}
\begin{picture}(5,4)(-5,-3.9)
  \put(0,0){\line(1,0){3}}
  \put(0,-1){\line(1,0){3}}
  \put(0,-2){\line(1,0){2}}
  \put(0,-3){\line(1,0){2}}
  \put(0,-4){\line(1,0){1}}
  \put(0,0){\line(0,-1){4}}
  \put(1,0){\line(0,-1){4}}
  \put(2,0){\line(0,-1){3}}
  \put(3,0){\line(0,-1){1}}
  \put(0.3,-.9){\small 6}
  \put(1.3,-.9){\small 4}
  \put(2.3,-.9){\small 1}
  \put(0.3,-1.9){\small 4}
  \put(1.3,-1.9){\small 2}
  \put(0.3,-2.9){\small 3}
  \put(1.3,-2.9){\small 1}
  \put(0.3,-3.9){\small 1}
\end{picture}
\end{center}
Given a partition $\lambda$ and a number $d$ we define the \emph{hook length} of $\lambda$
and the $d$-\emph{content} of $\lambda$ as \Mcambio the following polynomials:\finMcambio
\[
\hook_\lambda(q)=\prod_{u\in Y(\lambda)}[h(u)]_q \qquad
\Cnumbers^d_\lambda(q)=\prod_{u\in Y(\lambda)}[d+1+c(u)]_q.
\]

%
%

Let $\delta=(\delta_1,\delta_2)$ be a partition with at most two parts
and let $d=\delta_1-\delta_2$.
Let $\lambda$ be a partition with $\ell(\lambda)\le d+1$.
Since $s_\lambda$ is homogeneous of degree $|\lambda|$ it follows that  
\begin{equation*}
  s_\lambda(x_1^{\delta_1}x_2^{\delta_2},
x_1^{\delta_1-1}x_2^{\delta_2+1},...,
x_1^{\delta_2}x_2^{\delta_1})
=
(x_1^{\delta_1}x_2^{\delta_2})^{|\lambda|}s_\lambda(1,q,q^2,...,q^d), 
\end{equation*}
where $q=x_1^{-1}x_2$. 
If $b(\lambda)=\sum(i-1)\lambda_i$ and
\[
 P_\lambda^d(q)=\frac{\Cnumbers^d_\lambda(q)}{\hook_\lambda(q)}
\]
then {Theorem 7.21.2} in \cite{S} states that 
\begin{equation}\label{eq.S.7.21.2}
 s_\lambda(1,q,...,q^{d})=q^{b(\lambda)}P^d_\lambda(q).
\end{equation}
This identity is known as the \emph{Hook-content formula}, see
the notes in Ch. 7 of \cite{S} for more information about it.

It follows from \eqref{eq.Plethysm} that if  $x_1$ and $x_2$ are the  eigenvalues 
of $g\in\text{GL}(2,\CC)$, then
\begin{equation}\label{eq.CharPlethysm}
  \chi_{_{\SS_{\lambda}\left(\SS_{\delta}(\CC^2)\right)}}(g)
=(x_1^{\delta_1}x_2^{\delta_2})^{|\lambda|} q^{b(\lambda)}P^d_\lambda(q).
\end{equation}
%
%
%
%
%

%
%

\section{Main results}
\subsection{Equation \eqref{eq.Main} and the Hook-content formula}

The following theorem expresses the isomorphism condition 
of \eqref{eq.Main} in terms of the function $P$.
\begin{theorem}\label{Thm.SLGL}
 Let  $\delta=(\delta_1,\delta_2)$, $\epsilon=(\epsilon_1,\epsilon_2)$
and $d=\delta_1-\delta_2$, $e=\epsilon_1-\epsilon_2$.
Let $\lambda$, $\mu$ be partitions with $\ell(\lambda)\le d+1$ and $\ell(\mu)\le e+1$.
Then 
\medskip
\begin{enumerate}[(1)]
 \item $\SS_\lambda\left(\SS_\delta(\CC^2)\right)\simeq 
\SS_\mu\left(\SS_\epsilon(\CC^2)\right)$ as  $\text{SL}(2,\CC)$-modules if and only if 
\begin{equation}\label{eq.P}
P^d_\lambda=P^e_\mu
\end{equation}
and in this case $|\lambda|d-|\mu|e$ is even.
\medskip
\item  $\SS_\lambda\left(\SS_\delta(\CC^2)\right)\simeq 
\SS_\mu\left(\SS_\epsilon(\CC^2)\right)$ as  $\text{GL}(2,\CC)$-modules
if and only if in addition to \eqref{eq.P} 
it also holds
\begin{equation}\label{eq.N}
|\delta||\lambda|=|\epsilon||\mu|.
\end{equation}
\end{enumerate}
\end{theorem}

\medskip

\pf
On the one hand, it follows from \eqref{eq.CharPlethysm} that 
$
 \SS_{\lambda}\left(\SS_{\delta}(\CC^2)\right)\simeq \SS_{\mu}\left(\SS_{\epsilon}(\CC^2)\right)
$
 as representations of $\text{GL}(2,\CC)$ if and only if
\begin{equation}\label{eq.GL}
 (x_1^{\delta_1}x_2^{\delta_2})^{|\lambda|} q^{b(\lambda)}P^d_\lambda(q)=
(x_1^{\epsilon_1}x_2^{\epsilon_2})^{|\mu|} q^{b(\mu)}P^e_\mu(q)
\end{equation}
and since the identity $x_1x_2=1$ holds in $\text{SL}(2,\CC)$, it follows that 
$q=x_1^{-1}x_2=x_2^2$ and hence 
$
 \SS_{\lambda}\left(\SS_{\delta}(\CC^2)\right)\simeq \SS_{\mu}\left(\SS_{\epsilon}(\CC^2)\right)
$
as  $\text{SL}(2,\CC)$-modules if and only if
\begin{equation}\label{eq.SL}
x_2^{-d|\lambda|+2b(\lambda)}P^d_\lambda(x_2^{2})=
x_2^{-e|\mu|+2b(\mu)}P^e_\mu(x_2^{2})
\end{equation}
as a function of $x_2$.

On the other hand, since $s_\lambda$ is symmetric, it follows from 
\eqref{eq.Plethysm} and \eqref{eq.CharPlethysm} that 
\[
x_1^{\delta_1|\lambda|-b(\lambda)}x_2^{\delta_2|\lambda|+b(\lambda)}P^d_\lambda(q)
=
x_2^{\delta_1|\lambda|-b(\lambda)}x_1^{\delta_2|\lambda|+b(\lambda)}P^d_\lambda(q^{-1})
\]
and thus 
\begin{align*}
\frac{P^d_\lambda(q)}{P^d_\lambda(q^{-1})}
&=
x_2^{(\delta_1-\delta_2)|\lambda|-2b(\lambda)}x_1^{(\delta_2-\delta_1)|\lambda|+2b(\lambda)} \\ 
&=
q^{d|\lambda|-2b(\lambda)}.
\end{align*}
A similar identity holds for $\mu$ and $\epsilon$ instead of $\lambda$ and $\delta$.

We now assume condition \eqref{eq.P}. This and the above identities imply that 
\begin{equation}\label{eq.**}
 d|\lambda|-2b(\lambda)= e|\mu|-2b(\mu)
\end{equation}
and therefore \eqref{eq.SL} holds and thus  
$
 \SS_{\lambda}\left(\SS_{\delta}(\CC^2)\right)\simeq \SS_{\mu}\left(\SS_{\epsilon}(\CC^2)\right)
$
as representations of $\text{SL}(2,\CC)$. 
It also follows from \eqref{eq.**} that $|\lambda|d-|\mu|e$ is even.

If we additionally assume that condition \eqref{eq.N} holds, then 
adding and substracting  \eqref{eq.**} and \eqref{eq.N} we obtain
\begin{align*}
\delta_1|\lambda|-b(\lambda)&=\epsilon_1|\tau|-b(\tau) \\
\delta_2|\lambda|+b(\lambda)&=\epsilon_2|\tau|+b(\tau),
\end{align*}
and taking into account that $q=x_1^{-1}x_2$, \eqref{eq.GL} follows and thus 
$
 \SS_{\lambda}\left(\SS_{\delta}(\CC^2)\right)\simeq \SS_{\mu}\left(\SS_{\epsilon}(\CC^2)\right)
$
as representations of $\text{GL}(2,\CC)$. 

For the converse statements, we first observe that $q=0$ is neither a root nor a pole of the rational function $P_\lambda^d(q)=\frac{\Cnumbers^d_\lambda(q)}{\hook_\lambda(q)}$. 
Therefore, if 
$
 \SS_{\lambda}\left(\SS_{\delta}(\CC^2)\right)\simeq \SS_{\mu}\left(\SS_{\epsilon}(\CC^2)\right)
$
as representations of $\text{SL}(2,\CC)$ then it follows from \eqref{eq.SL} that 
$P^d_\lambda=P^e_\mu$. 
If the isomorphism also holds  
as representations of $\text{GL}(2,\CC)$, then we obtain \eqref{eq.N} by
specializing \eqref{eq.GL} at $x_1=x_2$. 
\eop


\subsection{GL$(2,\CC)$-isomorphisms from SL$(2,\CC)$-isomorphisms}\label{GfromS} 
Let  $\delta=(\delta_1,\delta_2)$, $d=\delta_1-\delta_2$, and let
$\lambda$ be partition with $\ell(\lambda)\le d+1$.    Since $d+1=\dim(\SS_{(d)}(\CC^2))$
it follows from the discusion in  \S\ref{subsec.PolyRep}, that if
$$
 \tilde\lambda=(\lambda_1-\lambda_{d+1},\dots,\lambda_{d}-\lambda_{d+1})
$$
then $\SS_\lambda(\left(\SS_{(d)}(\CC^2)\right))\simeq\SS_{\tilde\lambda}(\left(\SS_{(d)}(\CC^2)\right))$ as 
SL$(2,\CC)$-modules.  
Thus, in order to study the plethysm equation \eqref{eq.Main}
as SL$(V)$-modules it is enough to consider the problem of finding 
$d$, $e$, $\lambda$ and $\mu$ with 
 $\ell(\lambda)\le d$, $\ell(\mu)\le e$, such that 
\begin{equation}\label{eq.MainSL}
  \SS_\lambda\left(\SS_{(d)}(\CC^2)\right)\simeq \SS_\mu\left(\SS_{(e)}(\CC^2)\right)
\end{equation}
as representations of $\text{SL}(2,\CC)$.

On the other hand, if \eqref{eq.MainSL} holds,  
part (2) of {Theorem \ref{Thm.SLGL}} says that 
the isomorphism also holds as $\text{GL}(2,\CC)$-modules
if and only if $|\lambda|d=|\mu|e$.

If this is not the case,
a natural  question to ask is whether there exist $l,m,x,y\in\mathbb Z_{\ge 0}$ such that
\[
 \SS_{\lambda+(l^{d+1})}\left(\SS_{(d+x,x)}(\CC^2)\right)\simeq 
 \SS_{\mu+(m^{e+1})}\left(\SS_{(e+y,y)}(\CC^2)\right)
\]
as representations of $\text{GL}(2,\CC)$.

According to  part (2) of {Theorem \ref{Thm.SLGL}} the answer is positive 
if and only if 
\[
  (|\lambda|+l(d+1))(d+2x)=(|\mu|+m(e+1))(e+2y)
\]
wich is equivalent to 
\begin{equation}\label{eq.CondGL}
 \big(|\mu|+m(e+1)\big)y-\big(|\lambda|+l(d+1)\big)x 
 = \frac{|\lambda|d-|\mu|e}2+l\binom{d\!+\!1}2-m\binom{e\!+\!1}2.
\end{equation}
From part (1) of {Theorem \ref{Thm.SLGL}} we know that the 
right hand side of \eqref{eq.CondGL} is an integer number. 
In addition, there exist $l,m,x,y\in\mathbb Z_{\ge 0}$ satisfying 
\eqref{eq.CondGL} if and only if the 
there exist $l,m\in\mathbb Z_{\ge 0}$ such that 
\begin{equation}\label{eq.CondGL2}
 \text{gcd}\left\{\big(|\mu|+m(e+1)\big),\big(|\lambda|+l(d+1)\big)\right\}\left|
 \frac{|\lambda|d-|\mu|e}2+l\binom{d\!+\!1}2-m\binom{e\!+\!1}2\right..
\end{equation}
  
Such $l$ and $m$ do not always exist but in many cases they do.
Concretely

\begin{theorem} If $n$ is an integer number, let $\nu_2(n)$ be the exponent of the highest power of the prime 2 that divides $n$. 

Then there exist $l$ and $m$ such that 
 \eqref{eq.CondGL2} 
holds unless $\nu_2(|\mu|)\ne \nu_2(|\lambda|)$ and  $0<{\rm min}\{\nu_2(|\mu|),\nu_2(|\lambda|)\}<{\rm min}\{\nu_2(e+1),\nu_2(d+1)\}$.
\end{theorem}
Since this is a side issue with respect to the main thrust of this paper, and the proof, while not difficult, is slightly complicated, we will prove the above theorem in another article.

%
%
%

%
%

\subsection{Equation \eqref{eq.Main} as SL$(2,\CC)$-modules}

\begin{notation}
In order to write the proofs easier, if we have a rectangular array of $q$ numbers:

$$\overbrace{i\begin{cases}
\begin{matrix}
[x+i+j-2]&[x+i+j-3]&...&[x+i-1]\\
[x+i+j-3]&[x+i+j-3]&...&[x+i-2]\\
\vdots&&...&\vdots\\
 [x+j-1]&[x+j-2]&...&[x]\\
\end{matrix}
\end{cases}
}^{j}$$
\end{notation}

\noindent in which all the columns and rows  decrease by 1, 
we will denote  the product of all the elements $[*]$ in that rectangle by $\rho_{i,j}(x)$. 
 Clearly $\rho_{i,j}(x)=\rho_{j,i}(x)$ and 
if $k>j$ then $\rho_{i,k}(x)=\rho_{i,k-j}(x+j)\rho_{i,j}(x)$.

\begin{lemma}\label{hdetranspuesta}
If $\lambda^t$ is the transpose of $\lambda$ (see \S\ref{partitions}), then:
$$\hook_\lambda=\hook_{\lambda^t}$$
\end{lemma}
\pf
Let $x_1,\ldots ,x_t,y_1,...,y_t$ be such that the Young diagram of $\lambda $ is:
\begin{center}
 \setlength{\unitlength}{16pt}
\begin{picture}(5,5.0)(0,-4.1)
  \put(-3.4,-1){$Y(\lambda)=$}
  \put(0,0){\line(1,0){2}}
  \multiput(2,0)(.4,0){6}{\line(1,0){.2}}
  \put(4,0){\line(1,0){1}}
  \multiput(0,-1)(.4,0){10}{\line(1,0){.2}}
  \put(4,-1){\line(1,0){1}}
  \multiput(0,-2)(.4,0){10}{\line(1,0){.2}}
  \multiput(0,-3)(.4,0){5}{\line(1,0){.2}}
  \put(1,-3){\line(1,0){1}}
  \put(0,-4){\line(1,0){1}}
  \put(0,0){\line(0,-1){2}}
  \put(0,-3){\line(0,-1){1}}
  \multiput(0,-2)(0,-.3){5}{\line(0,-1){.12}}
  \multiput(1,0)(0,-.3){11}{\line(0,-1){.12}}
  \multiput(2,0)(0,-.3){10}{\line(0,-1){.12}}
  \multiput(4,0)(0,-.3){3}{\line(0,-1){.12}}
  \put(1,-3){\line(0,-1){1}}
  \put(4,-1){\line(0,-1){1}}
  \put(5,0){\line(0,-1){1}}
  \put(-0.8,-0.7){$y_1$}
  \put(-0.8,-1.7){$y_2$}
  \put(-0.8,-3.7){$y_{t}$}
  \put(0.2,.2){$x_1$}
  \put(1.2,.2){$x_2$}
  \put(4.2,.2){$x_t$}
\end{picture}
\end{center}
Then $\hook_\lambda$ is the product of all the $\rho_{y_i,x_j}(1+y_{i+1}+\dots +y_{t}+x_{j+1}+\dots +x_{t})$.
Since that product is obviously symmetric on the $x$'s and $y$'s, then we obtain the result.
\eop
\fincambio

\newcommand\fila [2]{#1_1,#1_2,\dots,#1_{#2}}
\newcommand\usual [3]{\left( (#1_1+#1_2+\cdots +#1_{#3-1}+#1_{#3})^{#2_1},(#1_1+#1_2+\cdots +#1_{#3-1})^{#2_2},\dots, (#1_1+#1_2)^{#2_{#3-1}},#1_1^{#2_{#3}}\right) }

\begin{notation}\label{minotacion}
Let $h_1,...,h_t,v_1,...,v_{t+1}$ be positive integers.
The { notation} $\dntg h1tv1t{t+1}$
 { will mean} 
the SL$(2,\CC)$-module
$\SS_\lambda\left(\SS_{(w)}(\CC^2)\right)$
where
$w=v_1+\cdots v_{r+1}-1$ and
\[
 \lambda=\big( (h_1+\cdots+h_{t-1} +h_{t})^{v_1},(h_1+\cdots +h_{t-1})^{v_2},
 \dots, (h_1+h_2)^{v_{t-1}},h_1^{v_{t}}\big),
 \]

In order to simplify this notation, given a sequence $x_1,x_2,...,x_t$, we will denote
by $\rvec x$ the sequence $x_1,x_2,...,x_t$ and by
$\lvec x$ the sequence $x_t,x_{t-1},...,x_1$.
That is:
\begin{center}
 \setlength{\unitlength}{16pt}
\begin{picture}(5,5.0)(0,-4.1)
\put(-4.8,-2.2){$\langle\rvec h||\rvec v\rangle=$}
  \put(0,0){\line(1,0){2}}
  \multiput(2,0)(.4,0){6}{\line(1,0){.2}}
  \put(4,0){\line(1,0){1}}
  \multiput(0,-1)(.4,0){10}{\line(1,0){.2}}
  \put(4,-1){\line(1,0){1}}
  \multiput(0,-2)(.4,0){10}{\line(1,0){.2}}
  \multiput(0,-3)(.4,0){5}{\line(1,0){.2}}
  \put(1,-3){\line(1,0){1}}
  \put(0,-4){\line(1,0){1}}
  \put(0,0){\line(0,-1){2}}
  \put(0,-3){\line(0,-1){1}}
  \multiput(0,-2)(0,-.3){5}{\line(0,-1){.12}}
  \multiput(1,0)(0,-.3){11}{\line(0,-1){.12}}
  \multiput(2,0)(0,-.3){10}{\line(0,-1){.12}}
  \multiput(4,0)(0,-.3){3}{\line(0,-1){.12}}
  \put(1,-3){\line(0,-1){1}}
  \put(4,-1){\line(0,-1){1}}
  \put(5,0){\line(0,-1){1}}
  \put(-0.8,-0.7){$v_1$}
  \put(-0.8,-1.7){$v_2$}
  \put(-0.8,-3.7){$v_r$}
  \put(0.2,.2){$h_1$}
  \put(1.2,.2){$h_2$}
  \put(4.2,.2){$h_r$}
\put(2,-4.3){$v_1+\dots+v_r+v_{r+1}$}
\end{picture}
\end{center}
\end{notation}
\fincambio

If 
In this notation,  
the  SL$(2,\CC)$-modules isomorphism given in Theorem \ref{eq.SL-Dual} becomes
\begin{theorem}\label{thm.complemento}
$$\langle\rvec h||\rvec v\rangle\simeq \langle\lvec h||\lvec v\rangle$$
\end{theorem}
In pictures:

 \setlength{\unitlength}{16pt}
\begin{picture}(5,5.0)(-2,-4.1)
  \put(0,0){\line(1,0){2}}
  \multiput(2,0)(.4,0){6}{\line(1,0){.2}}
  \put(4,0){\line(1,0){1}}
  \multiput(0,-1)(.4,0){10}{\line(1,0){.2}}
  \put(4,-1){\line(1,0){1}}
  \multiput(0,-2)(.4,0){10}{\line(1,0){.2}}
  \multiput(0,-3)(.4,0){5}{\line(1,0){.2}}
  \put(1,-3){\line(1,0){1}}
  \put(0,-4){\line(1,0){1}}
  \put(0,0){\line(0,-1){2}}
  \put(0,-3){\line(0,-1){1}}
  \multiput(0,-2)(0,-.3){5}{\line(0,-1){.12}}
  \multiput(1,0)(0,-.3){11}{\line(0,-1){.12}}
  \multiput(2,0)(0,-.3){10}{\line(0,-1){.12}}
  \multiput(4,0)(0,-.3){3}{\line(0,-1){.12}}
  \put(1,-3){\line(0,-1){1}}
  \put(4,-1){\line(0,-1){1}}
  \put(5,0){\line(0,-1){1}}
  \put(-0.8,-0.7){$v_1$}
  \put(-0.8,-1.7){$v_2$}
  \put(-0.8,-3.7){$v_r$}
  \put(0.2,.2){$h_1$}
  \put(1.2,.2){$h_2$}
  \put(4.2,.2){$h_r$}
\put(1.5,-4.3){$v_1+\dots+v_r+v_{r+1}$}
\put(7,-2){$\simeq$}
\end{picture}
\begin{picture}(5,5.0)(-6,-4.1)
  \put(0,0){\line(1,0){2}}
  \multiput(2,0)(.4,0){6}{\line(1,0){.2}}
  \put(4,0){\line(1,0){1}}
  \multiput(0,-1)(.4,0){10}{\line(1,0){.2}}
  \put(4,-1){\line(1,0){1}}
  \multiput(0,-2)(.4,0){10}{\line(1,0){.2}}
  \multiput(0,-3)(.4,0){5}{\line(1,0){.2}}
  \put(1,-3){\line(1,0){1}}
  \put(0,-4){\line(1,0){1}}
  \put(0,0){\line(0,-1){2}}
  \put(0,-3){\line(0,-1){1}}
  \multiput(0,-2)(0,-.3){5}{\line(0,-1){.12}}
  \multiput(1,0)(0,-.3){11}{\line(0,-1){.12}}
  \multiput(2,0)(0,-.3){10}{\line(0,-1){.12}}
  \multiput(4,0)(0,-.3){3}{\line(0,-1){.12}}
  \put(1,-3){\line(0,-1){1}}
  \put(4,-1){\line(0,-1){1}}
  \put(5,0){\line(0,-1){1}}
  \put(-1.5,-0.7){$v_{r+1}$}
  \put(-0.8,-1.7){$v_r$}
  \put(-0.8,-3.7){$v_2$}
  \put(0.1,.2){$h_r$}
  \put(1.0,.2){$h_{r-1}$}
  \put(4.2,.2){$h_1$}
\put(1.5,-4.3){$v_1+\dots+v_r+v_{r+1}$}
\end{picture}

\newcommand\suma [2]{#1_1+#1_2+\cdots+#1_{#2}}
\newcommand\reversefila [2]{#1_{#2},#1_{#2-1}\dots,#1_2,#1_1}
\newcommand\reversefilawoutlast [2]{#1_{#2},#1_{#2-1}\dots,#1_2}
%


%
%


\begin{theorem}\label{Main}
Let  $s\ge 0$ and $t=s$ or $t=s+1$.
Let $x_1,\dots,x_s$ and $y_1,\dots,y_t$ be two sequences of
positive integers, $u,v,z$ three  positive integers. Let $|\rvec x|$ denote $\sum_i x_i$.

a)
 The following
$\text{SL}(2,\CC)$-isomorphisms hold:
$$
\begin{matrix} 
\dnt {\rvec x}u{\rvec y}z{\rvec x}v{\rvec y}&\simeq& \dnt {\rvec x}v{\rvec y}z{\rvec x}u{\rvec y}\\[2mm]
\rule[-.36mm]{.6pt}{2.7mm}\,\wr && \rule[-.36mm]{.6pt}{2.7mm}\,\wr\\[2mm]
\dnt {\lvec y}u{\lvec x}{\lvec y}v{\lvec x}z&\simeq& \dnt {\lvec y}v{\lvec x}{\lvec y}u{\lvec x}z\\[2mm]
\end{matrix}
$$
b) 
Let $S=|\rvec x|^2+2\sum_{i,j: i+j=t}x_iy_j$. 
If $z(z-1)=S+|\rvec x|(u+v)$ in the case $t=s$, or
$z(z-1)=S+|\rvec x|(u+v)+uv$ in the case $t=s+1$
then the first row is a $\text{GL}(2,\CC)$-isomorphism.

c) Except in the trivial case $u=v$, the second row and both columns are never GL$(2,\CC)$ isomorphisms.
\end{theorem}
\pf
a) The horizontal isomorphisms reveal a symmetry between $u$ and $v$.  Since the vertical isomorphisms follow from Theorem \ref{thm.complemento}, we only need to prove one of the horizontal ones. We will prove the second one, i.e.,
we will show that $\dnt {\lvec y}u{\lvec x}{\lvec y}v{\lvec x}z$ is symmetric on $u$ and $v$.

 Let us call $\lambda_{u,v}$  the subyacent partition in $\dnt {\lvec y}u{\lvec x}{\lvec y}v{\lvec x}z$
Since $\lambda_{v,u}=\lambda_{u,v}^t$, then by Lemma \ref{hdetranspuesta} we have $\hook_{\lambda_{u,v}}=\hook_{\lambda_{v,u}}$.

Now let us see $\Cnumbers$.

We need to compute 
$\Cnumbers_{\lambda_{u,v}}^w$, where $w=|\lvec x|+v+|\lvec y|+z-1$. Note that $w$ depends on $v$ but not $u$.

In this case the product $\rro ijk.$ arises from an array  of the form:

$$\overbrace{i\begin{cases}
\begin{matrix}
[k+i-1]&[k+i]&...&[k+i+j-2]\\
\vdots&&...&\vdots\\
 [k+1]&[k+2]&...&[k+j]\\
 [k]&[k+1]&...&[k+j-1]\\
\end{matrix}
\end{cases}
}^{j}$$
Let's consider first the case $t=s$.
The partition is then: 
  
\begin{center}
 \setlength{\unitlength}{16pt}
\begin{picture}(5,8.0)(2,-7.1)
  \put(-4.3,-3){$Y(\lambda_{u,v}):$}
  \put(0,0){\line(1,0){7}}
  \put(0,-1){\line(1,0){7}}
  \put(0,-2){\line(1,0){6}}
  \put(0,-3){\line(1,0){5}}
  \put(0,-4){\line(1,0){4}}
  \put(0,-5){\line(1,0){3}}
  \put(0,-6){\line(1,0){2}}
  \put(0,-7){\line(1,0){1}}
  \put(0,0){\line(0,-1){7}}
  \put(1,0){\line(0,-1){7}}
  \put(2,0){\line(0,-1){6}}
  \put(3,0){\line(0,-1){5}}
  \put(4,0){\line(0,-1){4}}
  \put(5,0){\line(0,-1){3}}
  \put(6,0){\line(0,-1){2}}
  \put(7,0){\line(0,-1){1}}
  \put(0.2,.2){$y_s$}
  \put(1.2,.2){$\dots$}
  \put(2.2,.2){$y_1$}
  \put(3.2,.2){$u$}
  \put(4.2,.2){$x_s$}
  \put(5.2,.2){$\dots$}
  \put(6.2,.2){$x_1$}
  \put(-0.8,-0.7){$y_s$}
  \put(-0.8,-1.7){$\vdots$}
  \put(-0.8,-2.7){$y_1$}
  \put(-0.8,-3.7){$v$}
  \put(-0.8,-4.7){$x_s$}
  \put(-0.8,-5.7){$\vdots$}
  \put(-0.8,-6.7){$x_1$}
\end{picture}
\end{center}

We see from $Y(\lambda_{u,v})$ that $\Cnumbers_{\lambda_{u,v}}^w$ is the product of
\begin{enumerate}
\item    $\rho_{|\lvec y|+v,|\lvec y|+u}(w+2-|\lvec y|-v)$.
Note that $w+2-|\lvec y|-v=|\lvec x|+z+1$, this item
is $\rho_{|\lvec y|+v,|\lvec y|+u}(|\lvec x|+z+1)$,   thus symmetric in $u,v$.
\item $\rho$'s from the part of the table below the horizontal $v$ line, which are independent of $u,v$.
\item $\rho$'s from the part of the table to the right of vertical $u$ column.
These are of the form $\rho_{y_i,x_j}(*)$ and $*$ is of the form
   $w+ 2+|\lvec y|+u+$ some $x$'s $-$ some  $y$'s, i.e. $w+u+$ other stuff.
Note that since $w$ depends on $v$ and not $u$, then  $w+u$ is symmetric on $u,v$.
\end{enumerate}
Note that in the case $s=t=0$, the proof reduces to just the case $(1)$.

Now consider the case $t=s+1$. Now the partition is: 

\begin{center}
 \setlength{\unitlength}{16pt}
\begin{picture}(5,9.0)(2,-8.1)
  \put(-3.7,-3){$Y(\lambda_{u,v}):$}
  \put(0,0){\line(1,0){8}}
  \put(0,-1){\line(1,0){8}}
  \put(0,-2){\line(1,0){7}}
  \put(0,-3){\line(1,0){6}}
  \put(0,-4){\line(1,0){5}}
  \put(0,-5){\line(1,0){4}}
  \put(0,-6){\line(1,0){3}}
  \put(0,-7){\line(1,0){2}}
  \put(0,-8){\line(1,0){1}}
  \put(0,0){\line(0,-1){8}}
  \put(1,0){\line(0,-1){8}}
  \put(2,0){\line(0,-1){7}}
  \put(3,0){\line(0,-1){6}}
  \put(4,0){\line(0,-1){5}}
  \put(5,0){\line(0,-1){4}}
  \put(6,0){\line(0,-1){3}}
  \put(7,0){\line(0,-1){2}}
  \put(8,0){\line(0,-1){1}}
  \put(0.2,.2){$y_t$}
  \put(1.7,.2){$\dots$}
  \put(3.2,.2){$y_1$}
  \put(4.2,.2){$u$}
  \put(5.2,.2){$x_s$}
  \put(6.2,.2){$\dots$}
  \put(7.2,.2){$x_1$}
  \put(-0.8,-0.7){$y_t$}
  \put(-0.8,-2.2){$\vdots$}
  \put(-0.8,-3.7){$y_1$}
  \put(-0.8,-4.7){$v$}
  \put(-0.8,-5.7){$x_s$}
  \put(-0.8,-6.7){$\vdots$}
  \put(-0.8,-7.7){$x_1$}
\end{picture}
\end{center}

As in the previous case, the $\rho$'s from the part of the table below the horizontal $v$ line are independent of $u,v$
and the $\rho$'s from the part of the table to the right of vertical $u$ column depend on $u+v$ and thus are symmetric
on $u,v$. So the only problem is the central part of the table, which, unlike the previous case is not a rectangle.
Let's call $\Gamma$ the product of the $\rho$'s corresponding to $\Cnumbers_{\lambda_{u,v}}$ of that part of the table.
Here we simply observe that if we were to append an extra rectangle of height $v$ and length $u$
to the southeast corner of that part, then we would have a rectangle whose $\rho$, as in the previous case,
would be symmetric in $u,v$. But the rectangle appended contributes with a $\rho_{v,u}$ of something that does not depend on $u,v$, hence it is symmetric on $u,v$. Therefore $\Gamma$ is the quotient between two things symmetric on $u,v$, hence, symmetric itself.

b) Let $\mu_{u,v}$ now denote the subyacent partition in
$\dnt {\rvec x}u{\rvec y}z{\rvec x}v{\rvec y}$
and set $w_v=|\rvec x|+v+|\rvec y|+z-1$.
By part a) of this theorem and part b) of Theorem \ref{Thm.SLGL},
in order to prove $\dnt {\rvec x}u{\rvec y}z{\rvec x}v{\rvec y}\simeq\dnt {\rvec x}v{\rvec y}z{\rvec x}u{\rvec y}$ as  GL$(2,\CC)$-modules
it suffices to see that $|\mu_{u,v}|w_v=|\mu_{v,u}|w_u$, i.e. it is enough to see that $|\mu_{u,v}|w_v$ is symmetric in $u,v$.

Let us start with the case $s=t$.
The partition in this case is:

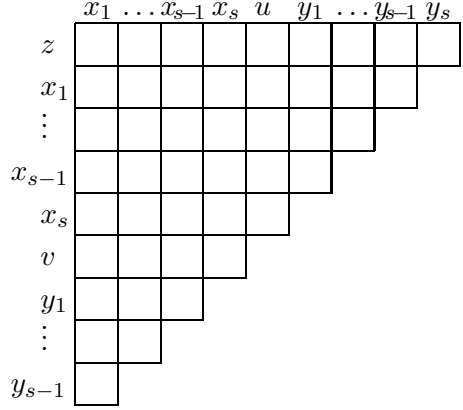
\begin{center}
 \setlength{\unitlength}{16pt}
\begin{picture}(5,9.0)(2,-8.1)
  \put(0,0){\line(1,0){9}}
  \put(0,-1){\line(1,0){9}}
  \put(0,-2){\line(1,0){8}}
  \put(0,-3){\line(1,0){7}}
  \put(0,-4){\line(1,0){6}}
  \put(0,-5){\line(1,0){5}}
  \put(0,-6){\line(1,0){4}}
  \put(0,-7){\line(1,0){3}}
  \put(0,-8){\line(1,0){2}}
  \put(0,-9){\line(1,0){1}}
  \put(0,0){\line(0,-1){9}}
  \put(1,0){\line(0,-1){9}}
  \put(2,0){\line(0,-1){8}}
  \put(3,0){\line(0,-1){7}}
  \put(4,0){\line(0,-1){6}}
  \put(5,0){\line(0,-1){5}}
  \put(6,0){\line(0,-1){4}}
  \put(7,0){\line(0,-1){3}}
  \put(8,0){\line(0,-1){2}}
  \put(9,0){\line(0,-1){1}}
  \put(0.2,.2){$x_1$}
  \put(1.1,.2){$\dots$}
 \put(2,.2){$x_{\!s\!-\!1\!}$}
  \put(3.2,.2){$x_s$}
  \put(4.2,.2){$u$}
  \put(5.2,.2){$y_1$}
  \put(6.1,.2){$\dots$}
  \put(7,.2){$y_{\!s\!-\!1\!}$}
  \put(8.2,.2){$y_s$}
  \put(-0.8,-0.7){$z$}
  \put(-0.8,-1.7){$x_1$}
  \put(-0.8,-2.7){$\vdots$}
  \put(-1.5,-3.7){$x_{s-1}$}
  \put(-0.8,-4.7){$x_s$}
  \put(-0.8,-5.7){$v$}
  \put(-0.8,-6.7){$y_1$}
  \put(-0.8,-7.7){$\vdots$}
  \put(-1.5,-8.7){$y_{s-1}$}
\end{picture}
\end{center}

\medskip

Thus $|\mu_{u,v}|=z(|\rvec x|+u+|\rvec y|)+|\rvec x|^2+|\rvec x|(u+v)+2\sum_{i,j: i+j=s}x_iy_j$.
Since $S=|\rvec x|^2+2\sum_{i,j: i+j=s}x_iy_j$.
(because this is the case $t=s$).
Hence

\begin{align*}
|\mu_{u,v}|w_v=\ &\left((|\rvec x|+u+|\rvec y|)z+S+|\rvec x|(u+v)\right).\left(|\rvec x|+v+|\rvec y|+z-1\right)\\
=\ &
(|\rvec x|+u+|\rvec y|)z(|\rvec x|+v+|\rvec y|)+
(|\rvec x|+|\rvec y|)z(z-1)+
\\
&+uz(z-1)+Sv+S(|\rvec x|+|\rvec y|+z-1)+|\rvec x|(u+v)v+\\
&+|\rvec x|(u+v)(|\rvec x|+|\rvec y|+z-1)
\end{align*}
The first, second, fifth and last terms are symmetric on $u,v$.
The third, fourth and sixth, using $z(z-1)=S+|\rvec x|(u+v)$, are equal to:
$$uz(z-1)+Sv+|\rvec x|(u+v)v=
(u+v)(S+|\rvec x|(u+v))$$ which is symmetric in $u,v$.

Let us analyze now the case $t=s+1$.
The partition in this case is:

\begin{center}
 \setlength{\unitlength}{16pt}
\begin{picture}(5,9.0)(2,-8.1)
  \put(0,0){\line(1,0){8}}
  \put(0,-1){\line(1,0){8}}
  \put(0,-2){\line(1,0){7}}
  \put(0,-3){\line(1,0){6}}
  \put(0,-4){\line(1,0){5}}
  \put(0,-5){\line(1,0){4}}
  \put(0,-6){\line(1,0){3}}
  \put(0,-7){\line(1,0){2}}
  \put(0,-8){\line(1,0){1}}
  \put(0,0){\line(0,-1){8}}
  \put(1,0){\line(0,-1){8}}
  \put(2,0){\line(0,-1){7}}
  \put(3,0){\line(0,-1){6}}
  \put(4,0){\line(0,-1){5}}
  \put(5,0){\line(0,-1){4}}
  \put(6,0){\line(0,-1){3}}
  \put(7,0){\line(0,-1){2}}
  \put(8,0){\line(0,-1){1}}
  \put(0.2,.2){$x_1$}
  \put(1.2,.2){$\dots$}
  \put(2.2,.2){$x_s$}
  \put(3.2,.2){$u$}
  \put(4.2,.2){$y_1$}
  \put(5.2,.2){$\dots$}
 \put(6.2,.2){$y_s$}
 \put(7.2,.2){$y_{s+1}$}
  \put(-0.8,-0.7){$z$}
  \put(-0.8,-1.7){$x_1$}
  \put(-0.8,-2.7){$\vdots$}
  \put(-0.8,-3.7){$x_s$}
  \put(-0.8,-4.7){$v$}
  \put(-0.8,-5.7){$y_1$}
  \put(-0.8,-6.7){$\vdots$}
  \put(-0.8,-7.7){$y_{s}$}
\end{picture}
\end{center}

Thus in this case  $|\mu_{u,v}|=z(|\rvec x|+u+|\rvec y|)+|\rvec x|^2+
|\rvec x|(u+v)+2\sum_{i,j: i+j=s+1}x_iy_j+uv$
Since in this case $2\sum_{i,j: i+j=s+1}x_iy_j=2\sum_{i,j: i+j=t}x_iy_j
$ we have
$|\mu_{u,v}|=z(|\rvec x|+u+|\rvec y|)+S+|\rvec x|(u+v)+uv$ and:

\begin{align*}
|\mu_{u,v}|w_v=\ &\left((|\rvec x|+u+|\rvec y|)z+S+|\rvec x|(u+v)+uv\right).\left(|\rvec x|+v+|\rvec y|+z-1\right)\\
=\ &
(|\rvec x|+u+|\rvec y|)z(|\rvec x|+v+|\rvec y|)+
(|\rvec x|+|\rvec y|)z(z-1)+
\\
&+uz(z-1)+Sv+S(|\rvec x|+|\rvec y|+z-1)+|\rvec x|(u+v)v+\\
&+|\rvec x|(u+v)(|\rvec x|+|\rvec y|+z-1)+uv(|\rvec x|+|\rvec y|+z-1)+uv^2\\
\end{align*}
The first, second, fifth, seventh and eight terms are symmetric on $u,v$.
The third, fourth, sixth and last term, using $z(z-1)=S+|\rvec x|(u+v)+uv$, are equal to $(u+v)(S+|\rvec x|(u+v)+uv)$, symmetric.

c)
The previous second horizontal isomorphism never holds as a GL$(2,\CC)$ isomorphism. (except in the trivial case $u=v$).

This follows since $\lambda_{v,u}=\lambda_{u,v}^t$, hence
$\vert\lambda_{u,v}\vert=\vert\lambda_{v,u}\vert$ but $w_v\ne w_u$
hence 
$\vert\lambda_{u,v}\vert w_v\ne \vert\lambda_{v,u}\vert w_u$, so by Theorem \ref{Thm.SLGL}
the GL$(2,\CC)$ isomorphism does not hold.
\eop

\begin{remark}
Note that by the first part of Theorem \ref{Thm.SLGL}, $\vert\lambda_{u,v}\vert w_v- \vert\lambda_{v,u}\vert w_u$ must be even. Since that diference is $\vert\lambda_{u,v}\vert(v-u)$ then either $v-u$ is even or $\lambda$ is even. This can also be verified directly.
\end{remark}

\begin{remark}
Although in the statement and proof of Theorem \ref{Main} all the variables must be positive, 
let us see what happens if we set some of them equal to 0.
\begin{itemize}
\item If we set one of the variables $y_i$ or $x_i$ equal to zero,
what happens is that this gives rise to another configuration with all variables positive but both $s$ and $t$ decrease by 1. For example, if we set $x_1=0$, this eliminates an $x$ variable, decreasing  $s$ by 1
but ``joins" $y_{s-1}$ and $y_s$ to form a new variable with value 
$y_{s-1}+y_s$, thus decreasing $s$ by 1 too. Hence the total number of variables decrease by two. 

\item If we set the variable $z=0$, then we decrease the 
total number of variables by 1, and we switch from the $t=s$ case to the $t=s+1$ case and viceversa, but we go for example from the upper isomorphsim of case $s=t$ to the lower isomorphism for case
$t=s+1$.
\end{itemize}

Therefore we could state just one theorem, in the following form:

\begin{theorem}\label{thm.MainMain}
Let  $s\ge 0$.
Let $x_1,\dots,x_s$ and $y_1,\dots,y_s$ be two sequences of
nonnegative integers, $u,v,z$ three  nonnegative integers . 
Then the following
$\text{SL}(2,\CC)$-isomorphisms hold:
$$
\begin{matrix} 
\dnt {\rvec x}u{\rvec y}z{\rvec x}v{\rvec y}&\simeq& \dnt {\rvec x}v{\rvec y}z{\rvec x}u{\rvec y}\\[2mm]
\rule[-.36mm]{.6pt}{2.7mm}\,\wr && \rule[-.36mm]{.6pt}{2.7mm}\,\wr\\[2mm]
\dnt {\lvec y}u{\lvec x}{\lvec y}v{\lvec x}z&\simeq& \dnt {\lvec y}v{\lvec x}{\lvec y}u{\lvec x}z\\[2mm]
\end{matrix}
$$
\end{theorem}
\end{remark}

%
%

\section{Some Corollaries}\label{Corollaries}

Here we obtain corollaries of Theorem \ref{Main}.

\begin{remark}
Hermite's is a corollary of our theorem, since it is the case $s=t=0$,
with $z=1$ which implies that the condition of part b) of Theorem \ref{Main}
is satisfied, since $z(z-1)=0$ while $S=|\rvec x|=0$ too. Hence we obtain
the full statement of Hermite's, while from Manivel's result only the SL$(2,\CC)$ isomorphism can be deduced.
\end{remark}

\newcommand\isoline [2]{\vvv {{#1}^s{#2}^{s+1}}{{#1}^{s+1}u{#2}^s}&\simeq& \vvv {{#1}^su{#2}^s}{{#1}^{s+1}{#2}^{s+1}}&\simeq& \vvv {{#1}^{s+1}{#2}^s}{{#1}^su{#2}^{s+1}}\\[2mm] }

\newcommand\reverseisoline [2]{\vvv {{#1}^{s+1}{#2}^s}{{#1}^su{#2}^{s+1}}&\simeq& \vvv {{#1}^su{#2}^s}{{#1}^{s+1}{#2}^{s+1}}&\simeq&\vvv {{#1}^s{#2}^{s+1}}{{#1}^{s+1}u{#2}^s}\\[2mm] }

\newcommand\isolinetwo [2]{\vvv {{#1}^su{#2}^{s+1}}{{#1}^{s+2}{#2}^{s+1}}&\!\!\!\simeq\!\!\!& \vvv {{#1}^{s+1}{#2}^{s+1}}{{#1}^{s+1}u{#2}^{s+1}}&\!\!\!\simeq\!\!\!& \vvv {{#1}^{s+1}u{#2}^s}{{#1}^{s+1}{#2}^{s+2}}\\[2mm] }

\newcommand\reverseisolinetwo [2]{ \vvv {{#1}^{s+1}u{#2}^s}{{#1}^{s+1}{#2}^{s+2}}&\!\!\!\simeq\!\!\!& \vvv {{#1}^{s+1}{#2}^{s+1}}{{#1}^{s+1}u{#2}^{s+1}}&\!\!\!\simeq\!\!\!& \vvv {{#1}^su{#2}^{s+1}}{{#1}^{s+2}{#2}^{s+1}}\\[2mm] }

\begin{theorem}
Let $v,z,u$ be positive integers. Let $s\ge 0$. Then the two following families of isomorphism hold:
\begin{equation}\tag{I}
\begin{matrix}
\isoline zv
\rule[-.36mm]{.6pt}{2.7mm}\,\wr && \rule[-.36mm]{.6pt}{2.7mm}\,\wr&& \rule[-.36mm]{.6pt}{2.7mm}\,\wr\\[2mm]
\isoline vz
\end{matrix}
\end{equation}
and
\begin{equation}\tag{II}
\begin{matrix} 
\isolinetwo zv
\rule[-.36mm]{.6pt}{2.7mm}\,\wr && \rule[-.36mm]{.6pt}{2.7mm}\,\wr&& \rule[-.36mm]{.6pt}{2.7mm}\,\wr\\[2mm]
\isolinetwo vz
\end{matrix}
\end{equation}
\end{theorem}
\pf
From Theore \ref{Main} we obtain:
\[
\begin{matrix} 
\vvv {{z}^su{v}^s}{{z}^{s+1}{v}^{s+1}}&\simeq&\vvv {{z}^s{v}^{s+1}}{{z}^{s+1}u{v}^s}\\[2mm]
 \rule[-.36mm]{.6pt}{2.7mm}\,\wr && \rule[-.36mm]{.6pt}{2.7mm}\,\wr\\[2mm]
\vvv {{v}^su{z}^s}{{v}^{s+1}{z}^{s+1}}&\simeq& \vvv {{v}^{s+1}{z}^s}{{v}^su{z}^{s+1}}\\[2mm] 
\end{matrix}
\]

If we applly the lower isomorphism to $\vvv {{z}^su{v}^s}{{z}^{s+1}{v}^{s+1}}$ we obtain:

\[
\begin{matrix} 
\reverseisoline zv
&&\rule[-.36mm]{.6pt}{2.7mm}\,\wr && \rule[-.36mm]{.6pt}{2.7mm}\,\wr\\[2mm]
&& \vvv {{v}^su{z}^s}{{v}^{s+1}{z}^{s+1}}&\simeq& \vvv {{v}^{s+1}{z}^s}{{v}^su{z}^{s+1}}\\[2mm] 
\end{matrix}
\]
Applying theorem \ref{thm.complemento} to the top left, we get:

$$
\begin{matrix} 
\reverseisoline zv
\rule[-.36mm]{.6pt}{2.7mm}\,\wr && \rule[-.36mm]{.6pt}{2.7mm}\,\wr&& \rule[-.36mm]{.6pt}{2.7mm}\,\wr\\[2mm]
\isoline vz
\end{matrix}
$$
Rearranging the top line we get the first result.
The second result is similar: from Theorem \ref{Main} we get:

\[
\begin{matrix} 
\vvv {{z}^{s+1}{v}^{s+1}}{{z}^{s+1}u{v}^{s+1}}&\simeq& \vvv {{z}^su{v}^{s+1}}{{z}^{s+2}{v}^{s+1}}\\[2mm] 
\rule[-.36mm]{.6pt}{2.7mm}\,\wr && \rule[-.36mm]{.6pt}{2.7mm}\,\wr\\[2mm]
\vvv {{v}^{s+1}{z}^{s+1}}{{v}^{s+1}u{z}^{s+1}}&\simeq& \vvv {{v}^{s+1}u{z}^s}{{v}^{s+1}{z}^{s+2}}\\[2mm]
\end{matrix}
\]
Again, applying the lower isomorphism to $\vvv {{z}^{s+1}{v}^{s+1}}{{z}^{s+1}u{v}^{s+1}}$,
using theorem \ref{thm.complemento} and rearranging the top line we get the result.
\eop

\begin{remark}
Note that the isomorphism $I$ with $s=0$ gives:
$$
\begin{matrix} 
\vvv v{zu}&\simeq&\vvv u{zv}&\simeq& \vvv z{vu}\\[2mm] 
\rule[-.36mm]{.6pt}{2.7mm}\,\wr && \rule[-.36mm]{.6pt}{2.7mm}\,\wr&& \rule[-.36mm]{.6pt}{2.7mm}\,\wr\\[2mm]
\vvv z{vu}&\simeq&\vvv u{vz}&\simeq& \vvv v{zu}\\[2mm]
 \end{matrix}
$$
i.e., it says that $\vvv z{vu}$ is symmetric in $z,v,u$. This is Manivel's result.
It is not possible to say more because in order to apply \ref{Main} b) in this case,
we would need $z(z-1)=0$ which only happens when $z=1$, and this is
Hermite's.
\end{remark}

\begin{remark}
The isomorphism $II$ with $s=0$ gives:
$$
\begin{matrix} 
\vvv {uv}{zzv}&\simeq&\vvv {zv}{zuv}&\simeq& \vvv {zu}{zvv}\\[2mm] 
\rule[-.36mm]{.6pt}{2.7mm}\,\wr && \rule[-.36mm]{.6pt}{2.7mm}\,\wr&& \rule[-.36mm]{.6pt}{2.7mm}\,\wr\\[2mm]
\vvv {uz}{vvz}&\simeq&\vvv {vz}{vuz}&\simeq& \vvv {vu}{vzz}\\[2mm] 
\end{matrix}
$$
and the topright isomorphism is a GL$(2,\CC)$ isomorphism if $z(z-1)=uv$.
\end{remark}

\begin{remark}
If $s\ge 1$ we cannot obtain GL$(2,\CC)$ isomorphisms from the isomorphisms I or II.
\end{remark}


\end{document}